\documentclass{article}

\usepackage{PRIMEarxiv}

\usepackage[utf8]{inputenc} 
\usepackage[T1]{fontenc}    
\usepackage{hyperref}       
\usepackage{url}            
\usepackage{booktabs}       
\usepackage{amsfonts}       
\usepackage{amsmath}
\usepackage{nicefrac}       
\usepackage{microtype}      
\usepackage{lipsum}
\usepackage{fancyhdr}       
\usepackage{graphicx}       
\graphicspath{{media/}}     

\title{Existence of limit of multivariable function}

\author{
  Ming Yang\\
  \texttt{} \\
}

\begin{document}
\maketitle

\begin{abstract}
This article provides several theorems regarding the existence of limit for multivariable function, among which Theorem 1 and Theorem 3 relax the requirement for sequence of Heine's definition of limit. These results address the question of which paths need to be considered to determine the existence of limit for multivariable function.
\end{abstract}


\section{Introduction}
Unlike the limit of a single-variable function, the limit of a multivariable function appears to be more complex. A concise statement is that the limit of a multivariable function exists at a point only when the limit along any path to that point exists and is equal. However, the question of whether this condition can be relaxed has aroused interest. It is well-known that ensuring the existence of limits along all rays to a point does not guarantee the existence of the limit of a multivariable function at that point. In fact, it has been proven that even when the limits along any path of the form $y=cx^{m/n}$ or $x^m=(y/c)^n$ ($m$ and $n$ are relatively prime positive integers and $c$ is a nonzero constant) to the origin exist and are equal, it still does not ensure the existence of the limit of a binary function at the origin \cite{Real_Analysis1}. So, how many path limits are needed to guarantee the limit of a multivariable function?

\section{Theorem}
\label{sec:headings}

In this section, we will establish several theorems related to the aforementioned problem, thereby providing answers. Theorem 1 is a spacial version of Heine's definition of limit \cite{mahmudov2013single}, demonstrating the connection between the limit of binary function and sequence. Theorem 2 further indicates the relationship between the limit of binary function and limit along path. Theorem 3 extends Theorem 1 to the general case of multivariable function through non-constructive proof method.

\subsection{"Constructive" proof for simple case}
\noindent \textbf{Theorem 1.} For any constant $L$ and any accumulation point $P_0:(x_0,y_0)$ of $D\subset\mathbf{R}^2$, the necessary and sufficient condition for $\lim\limits_{\substack{P \to P_0 \\ P \in D}} f(P) = L$ is: for any constant $\phi_0 \in [0, 2\pi]$ and sequences $\{r_n\}$, $\{\phi_n\}$ such that $\lim\limits_{{n \to \infty}} r_n = 0$ and $\lim\limits_{{n \to \infty}} \phi_n = \phi_0$, the sequence of points $\{P_n : (x_0 + r_n \cos\phi_n, y_0 + r_n \sin\phi_n)\} \subset D$ satisfies $\lim\limits_{\substack{{n\to\infty}\\{P_n\in D}}} f(P_n) = L$.

\noindent \textbf{\textit{Proof.}} The necessity is easily proven, we only need prove its sufficiency. The definition of $\lim\limits_{\substack{P \to P_0 \\ P \in D}} f(P) = L$ states that $P_0$ is an accumulation point of $D$ and for any given $\varepsilon$, there exists a punctured neighborhood $U_0(P_0,r)$ such that for any point $P(x, y) \in U_0(P_0,r) \cap D$, we have $|(f(x, y) - l)| < \varepsilon$. Therefore, if $\lim\limits_{\substack{P \to P_0 \\ P \in D}} f(P)$ does not exist or is not equal to $L$, then $\exists \varepsilon$ is such that this neighborhood does not exist. So, take a point $P_1: (x_1, y_1)$ such that $| f(x_1, y_1) - l | \geq \varepsilon$, and within the neighborhood of radius $|( \vec{P_1} - \vec{P_0} )/2|$ of $P_0$, another point $P_2: (x_2, y_2)$ can be obtained such that $| f(x_2, y_2) - l | \geq \varepsilon$ (since $P_0$ is an accumulation point of $D$). Similarly, an infinite sequence $\{P_1, P_2, \ldots, P_n, \ldots\}$ can be constructed. The recursive relation implies $|(\vec{P_n} - \vec{P_0})| < |(\vec{P_1} - \vec{P_0})/2^{(n-1)}|$, so
\begin{equation}\label{eq:1}
\lim\limits_{{n \to \infty}} P_n = P_0
\end{equation}
Now, consider representing points on the plane as $(x_0 + r\cos\phi, y_0 + r\sin\phi)$, where $P_n$ can be expressed as $(x_0 + r_n\cos\phi_n, y_0 + r_n\sin\phi_n)$. This gives rise to an infinite sequence $\{\phi_n\}$. Let's define the domain of $\phi$ as $[0,2\pi]$, which can be divided into intervals $[0,\pi]$ and $[\pi,2\pi]$, i.e., $[0,2\pi]=[0,\pi] \cup [\pi,2\pi]$. We denote an interval containing infinitely many elements from $\{\phi_n\}$ as $I_1$, and choose a value $\phi_{i_1}$ from $\{\phi_n\}$ within $I_1$, corresponding to a point $P_{i_1}$ in ${P_n}$. Using bisection to further bisect intervals, for example, if $I_1$ is $[0,\pi]$, we divide it into $[0,\pi/2]$ and $[\pi/2,\pi]$, creating another interval $I_2$ containing infinitely many elements from $\{\phi_n\}$. This allows us to obtain a value $\phi_{i_2}$ in $I_2$ such that $i_2>i_1$, corresponding to a point $P_{i_2}$ in ${P_n}$. By repeating this process, an infinite sequence $\{I_1, I_2, \ldots, I_n, \ldots\}$ is constructed, ensuring that any interval $I_n$ contains infinitely many $\phi_n$, and we obtain an infinite sequence $\{\phi_{i_n}\}$ along with the corresponding points $\{P_{i_n}\}$. Assuming a constant parameter $\phi_0$ included in all closed intervals, let $\phi^{\prime}=\phi-\phi_0$. With the recursive relation $|I_{n+1}|=|I_n|/2$, we have $|I_n|=\pi/2^{n-1}$, and both $\phi_0$ and $\phi_{i_n}$ are contained in $I_n$, it follows that $-\pi/2^{n-1}\leq \phi^{\prime}_{i_n}=\phi_{i_n}-\phi_0\leq \pi/2^{n-1}$, leading to
\begin{equation}
\lim\limits_{{n\to\infty}}\phi_{i_n} = \lim\limits_{{n\to\infty}} \phi^{\prime}_{i_n}+\phi_0 = \phi_0
\end{equation}
Since integer $i_n>i_{n-1}$ implies $i_n\geq n$, we have $
0< r_{i_n}\leq r_n< r_1/2^{n-1}$, so
\begin{equation}
\lim\limits_{{n\to\infty}} r_{i_n}=0
\end{equation}
Therefore, if $\lim\limits_{\substack{P \to P_0 \\ P \in D}} f(P)$ does not exist or is not equal to $L$, then there always exists a sequence of points ${P_{i_n}}$ satisfying the condition $\phi_0\in[0,2\pi]$, $\lim\limits_{n\to\infty} r_{i_n} = 0$, $\lim\limits_{n\to\infty}\phi_{i_n} = \phi_0$, and $\lim\limits_{\substack{{n\to\infty}\\{P_{i_n}\in D}}} f(P_{i_n}) \neq L$ or not exist. Therefore, if any sequence of points satisfying the above conditions also satisfies $\lim\limits_{\substack{{n\to\infty}\\{P_n\in D}}} f(P_n) = L$, then $\lim\limits_{\substack{P \to P_0 \\ P \in D}} f(P) = L$. The sufficiency is proven, so Theorem 1 holds.

\noindent \textbf{Theorem 2.} For any constant $L$ and any accumulation point $P_0$ of $D\subset\mathbf{R}^2$, if there exists a punctured neighborhood $U_0(P_0,r)$ satisfy $(U_0(P_0,r)\cap D) \cup P_0$ is convex set, then the necessary and sufficient condition for the limit  $\lim\limits_{\substack{P \to P_0 \\ P \in D}} f(P) = L$ is: for any constant $\phi_0\in[0,2\pi]$ and function $\phi(r)$ satisfying $\lim\limits_{{r\to 0}} \phi(r) = \phi_0$, constructing the path $(P(r) : (x_0+r\cos\phi(r), y_0+r\sin\phi(r))) \subset D$, it holds that $\lim\limits_{\substack{r \to 0 \\ P(r) \in D}} f(P(r)) = L$.

\noindent \textbf{\textit{Proof.}} Similarly, the necessity is easily proven, it suffices to demonstrate sufficiency. If $\lim\limits_{\substack{P \to P_0 \\ P \in D}} f(P)$ does not exist or is not equal to $L$, since there exists an neighborhood $U_0(P_0,r)$ satisfy $(U_0(P_0,r)\cap D) \cup P_0$ is convex set, let's consider the construction method in Theorem 1. Take a sequence of points within this neighborhood, denoted as $\{P_{i_n}\}$. Consequently, we have
\begin{equation}
|\phi_{i_n} - \phi_{i_{n+1}}| \leq |I_n| = \frac{\pi}{2^{n-1}}
\end{equation}
So, note that when $n \geq 3$, $|\phi_{i_n} - \phi_{i_{n+1}}| \leq \frac{\pi}{4}$. Construct the triangle $P_0 P_{i_n} P_{i_{n+1}}$, which implies that the angle $\angle P_{i_n} P_0 P_{i_{n+1}} \leq \frac{\pi}{4}$, and the ratio of the sides
\begin{equation}
\frac{\left|\overline{P_0 P_{i_{n}}}\right|}{\left|\overline{P_0 P_{i_{n+1}}}\right|} = \frac{r_{i_{n}}}{r_{i_{n+1}}} > 2
\end{equation}
Let's denote the vertices of the triangle as $P_0, P_{i_n}, P_{i_{n+1}}$ with corresponding internal angles $A, B, C$, and the respective sides as $a, b, c$. Now $A \leq \frac{\pi}{4}$, and $c>2b$, then we have
\begin{equation}
a^2 = b^2 + c^2 - 2bc \cos A \leq b^2 + c^2 - \sqrt{2}bc
\end{equation}
\begin{equation}
\cos C = \frac{a^2 + b^2 - c^2}{2ab} \leq \frac{2b^2 - \sqrt{2}bc}{2ab} < \frac{(1 - \sqrt{2})b}{a} < 0
\end{equation}
Therefore, the angle $P_{i_n} P_{i_{n+1}} P_0$ is obtuse angle, indicating that the distance from a point on the line segment $P_{i_n} P_{i_{n+1}}$ (from $P_{i_n}$ to $P_{i_{n+1}}$) to $P_0$ is monotonically decreasing. Consequently, noticed \eqref{eq:1}, any point on the path $(P_{i_3}, P_{i_4}, \ldots, P_{i_n}, P_{i_{n+1}}, \ldots)$ can be represented as $P(r): (x_0 + r \cos\phi(r), y_0 + r \sin\phi(r))$. Since the angular deviation $\phi(r)$ of any point $P(r)$ on the line segment $P_{i_n} P_{i_{n+1}}$ and the angle $\phi_0$ are both contained within $I_n$, and if $m > n$, then $I_m$ is also contained within $I_n$, for any given $\varepsilon > 0$, there always exists $r_{i_n}$ such that when $r < r_{i_n}$, we have $-\varepsilon < -\frac{\pi}{2^{n-1}} \leq \phi^{\prime}(r) = \phi(r) - \phi_0 \leq \frac{\pi}{2^{n-1}} < \varepsilon$. Therefore,
\begin{equation}
\lim\limits_{{r \to 0}} \phi(r) = \phi_0
\end{equation}
Since $(U_0(P_0,r)\cap D) \cup P_0$ is convex, any point on the path belongs to $D$. Due to $\lim\limits_{\substack{{n\to\infty}\\{P_n\in D}}} f(P_{i_n}) \neq L$ or not exist, it is evident that $\lim\limits_{\substack{r \to 0 \\ P(r) \in D}} f(P(r))$ does not exist or is not equal to $L$. Therefore, if any path satisfying the given conditions also satisfies $\lim\limits_{\substack{r \to 0 \\ P(r) \in D}} f(P(r)) = L$, then $\lim\limits_{\substack{P \to P_0 \\ P \in D}} f(P) = L$. Sufficiency is thus proven, and Theorem 2 holds.

Theorem 2 requires the existence of an neighborhood of \( P_0 \) such that $(U_0(P_0,r)\cap D) \cup P_0$ forms a convex set. In fact, this condition can be relaxed to a more general connected set, as long as it is possible to construct path $P(s)$ ($s$ is the distance from the corresponding point to \( P_0 \) on this path) in $D$ connecting the sequence of points \( \{ P_{i_n} \} \) and corresponding function $\phi(s)$ satisfying $\lim\limits_{{s\to 0}} \phi(s) = \phi_0$ like Theorem 2. It meaning that the tangent to the path through \( P_0 \) has angle with x-axis of $\phi_0$. From a intuitive perspective, this implies that it is sufficient to consider the limits of the function along all paths passing through \( P_0 \) which exists a tangent line at \( P_0 \). This allows us to determine the limit of the function at the point. This serves as our answer to the problem in $\mathbf{R}^2$ case posed in the introduction.

\subsection{Non-constructive proof for general case}
Noticing the similarity between Theorem 1 and Bolzano–Weierstrass theorem \cite{Real_Analysis2}, we realize that Theorem 1 can be extended to a more general case and proven non-constructively through Bolzano–Weierstrass theorem. For any point $P\in\mathbf{R}^n$, it is always possible to find a set of coordinates $(x^1,x^2,\ldots,x^n)$ to represent that point and the coordinate set has a subset $\{x^{i_1},x^{i_2},\ldots,x^{i_{n'}}\}$ such that all elements are bounded in the punctured neighborhood of $P$. For instance, in hyperspherical coordinates centered at $P$ : $(r,\phi_1,\phi_2,\ldots,\phi_{n-1})$, any non-empty subset of $\{r,\phi_1,\phi_2,\ldots,\phi_{n-1}\}$ satisfies the given conditions. Then, under this set of coordinates and its subsets concerning the point $P_0$, the following theorem can be stated:

\noindent \textbf{Theorem 3.} For any constants $L$ and any accumulation point $P_0$ of $D\subset\mathbf{R}^n$, the necessary and sufficient condition for the function limit $\lim\limits_{\substack{P \to P_0 \\ P \in D}} f(P) = L$ is that: for any point sequence $\{P_m:(x^1,x^2,\ldots,x^n)\} \subset D$ such that all elements of $\{x_m^{i_1},x_m^{i_2},\ldots,x_m^{i_{n'}}\}$ converge and $\lim\limits_{{m \to \infty}} P_m = P_0$, it hold that $\lim\limits_{\substack{{m\to\infty}\\{P_m\in D}}} f(P_m) = L$.

\noindent \textbf{\textit{Proof.}} The necessity is easily proven; we now proceed to prove sufficiency. Following the method in Theorem 1, construct a sequence of points $\{P_m\}$ in the punctured neighborhood of $P_0$ satisfying the conditions. All elements of $\{x_m^{i_1}, x_m^{i_2}, \ldots, x_m^{i_{n'}}\}$ are bounded. By the Bolzano–Weierstrass theorem, the set of corresponding sequences $\{x_m^i\}$ contains convergent subsequences, meaning all elements of $\{x_m^{i_1}, x_m^{i_2}, \ldots, x_m^{i_{n'}}\}$ converge. Let $\{x_m^i\}$ have a convergent subsequence $\{x_{j_m}^i\}$. Construct the corresponding subsequence $\{P_{j_m}\}$ of $\{P_m\}$ such that all elements of $\{x_{j_m}^{i_1}, x_{j_m}^{i_2}, \ldots, x_{j_m}^{i_{n'}}\}$ converge. Similar to \eqref{eq:1}, $\lim\limits_{{m \to \infty}} P_m = P_0$, so $\lim\limits_{{m \to \infty}} P_{j_m} = P_0$. Therefore, if $\lim\limits_{\substack{P \to P_0 \\ P \in D}} f(P)$ does not exist or is not equal to $L$, there always exists a sequence $\{P_{j_m}\}$ satisfying the conditions, while $\lim\limits_{\substack{{m\to\infty}\\{P_m\in D}}} f(P_{j_m}) \neq L$ or not exist. Hence, if all satisfying sequences have $\lim\limits_{\substack{{m\to\infty}\\{P_m\in D}}} f(P_m) = L$, then $\lim\limits_{\substack{P \to P_0 \\ P \in D}} f(P) = L$. Sufficiency is proven, and Theorem 3 is established.

The relationship between Theorem 1 and Theorem 2 suggests that the feasibility of extending Theorem 2 to $\mathbf{R}^n$ case based on Theorem 3. In the example of hyperspherical coordinates, one can  establish the limit of multivariable function by the limit of the function along all paths passing through $P_0$ which exists a tangent line at $P_0$.

\section{Conclusion}
First, we prove Theorem 1, where a method for constructing non-convergent sequences satisfying the conditions is provided. It is worth noting that, in some cases, the construction method here is impractical. Indeed, considering Theorem 1 as a Heine's definition of limit under relaxed condition, it is evident that a truly constructive proof of Theorem 1 is impossible. Only by admitting a weak form of the axiom of choice can the equivalence of the Heine's definition of limit and the typical definition of limit be established \cite{moore2011early}, thus making it impossible to constructively prove the equivalence of the two definitions \cite{bridges1997constructive}. The sufficiency part of proving Theorem 1 can be seen as proving the sufficiency of the Heine's definition of limit, and the necessity of which is obviously valid. Therefore, if there exists a constructive proof of Theorem 1, the equivalence of the two limit definitions would also hold, which is not the case. Secondly, we prove Theorem 2, which essentially give the relationship between the limit of binary function and the limit of single variable function (though it is actually a special form of Theorem 5 in article \cite{darriba2021translated}). Finally, we prove Theorem 3, further extending Theorem 1 to the general case of multivariable function.

\section*{Acknowledgments}
The author proved Theorem 1 several years ago while studying calculus. With the development of AI translation tools like GPT, there is now an opportunity to fully document it in an English paper.

\bibliographystyle{unsrt}  
\bibliography{references}

\end{document}